\documentclass[11pt]{article}

\usepackage{graphicx}
\usepackage{times}
\usepackage{bm}
\usepackage{graphicx}
\usepackage{amsfonts}
\usepackage{amsmath}
\usepackage{amsthm}
\usepackage{amssymb}
\usepackage{natbib}
\usepackage{cite}
\usepackage{enumerate}
\usepackage[mathscr]{euscript}

\makeindex
\newcommand{\Pb}{\mathbb{P}}
\newcommand{\E}{\mathbb{E}}
\newcommand{\comment}[1]{}

\newtheorem{theorem}{Theorem}

\newtheorem{corollary}{Corollary}
\newtheorem{lemma}{Lemma}
\begin{document}

\title{AR(1) sequence with random coefficients: Regenerative properties and its application}
% Random Coefficient Autoregressive Process}

\author{
\hspace{.2in} \parbox[t]
{0.3\textwidth}
{{\sc Krishna B. Athreya}
\thanks{Departments of Mathematics and Statistics,
        Iowa State University, Iowa 50011, USA and Distinguished Visiting Profrssor, Department of Mathematics, IIT Bombay, Mumbai 400076, India. kba@iastate.edu.} }
%\hspace{.002\textwidth}
\parbox[t]{0.3\textwidth}
{{\sc Koushik Saha}
\thanks{Department of Mathematics,
        Indian Institute of Technology,
         Bombay, Mumbai 400076, India.  ksaha@math.iitb.ac.in.}}
%\hspace{.002\textwidth}
\parbox[t]{0.4\textwidth}{{\sc Radhendushka Srivastava}
\thanks{ Department of Mathematics,
        Indian Institute of Technology,
         Bombay, Mumbai 400076, India.  radhe@math.iitb.ac.in.}} }

\maketitle

\begin{abstract}
  Let $\{X_n\}_{n\ge0}$ be a sequence of real valued random variables such that $X_n=\rho_n X_{n-1}+\epsilon_n,~n=1,2,\ldots$, where $\{(\rho_n,\epsilon_n)\}_{n\ge1}$ are i.i.d. and independent of initial value (possibly random) $X_0$. In this paper it is shown that, under some natural conditions on the distribution of $(\rho_1,\epsilon_1)$, the sequence $\{X_n\}_{n\ge0}$ is regenerative in the sense that it could be broken up into i.i.d. components. Further, when $\rho_1$ and $\epsilon_1$ are independent,  we construct a non-parametric  strongly consistent estimator of the characteristic functions of $\rho_1$ and $\epsilon_1$.
\end{abstract}

\section{Introduction}
Let $\{X_n\}_{n\ge 0}$ be a sequence of real valued random variables satisfying the stochastic recurrence equation
\begin{equation}\label{model}
X_{n}=\rho_{n}X_{n-1}+\epsilon_{n},\quad n=1,2,\ldots
\end{equation}
where $\{(\rho_n,\epsilon_n)\}_{n\ge 1}$ are i.i.d.  $\mathbb{R}^2$-valued random vectors and independent of the initial random variable $X_0$. If
%$E|\epsilon_n|<\infty$,
$E(|X_0|)<\infty$ and $E(\epsilon_n)=0$, for each $n\ge 1$ and  then $E(X_n|X_0,\ldots,X_{n-1})=\rho_n X_{n-1}$. For this reason the sequence $\{X_n\}$ satisfying \eqref{model} is often referred to in the time series literature  as Random Coefficient Auto Regressive sequence of order one (RCAR(1))
%Several researchers have worked on RCAR models
(see \citet{A1976,R1978,NQ1980,B1986}). \citet{AHS2006} studied a parametric model  for $(\rho_1,\epsilon_1)$ under the assumption that $\rho_1$ and $\epsilon_1$ are independent and provided a consistent estimator of the model  parameters.
%Here we do not assume $\E(\epsilon_n)=0$. Instead we  focus on the {\it regenerative} property of the sequence $\{X_n\}$ satisfying \eqref{model}. We also study the convergence in distribution of $\{X_n\}$
In the current paper, we find conditions on the distribution function of $(\rho_1,\epsilon_1)$ to ensure that $\{X_n\}$  is a Harris recurrent Markov chain  and hence regenerative, i.e., it can be broken up into i.i.d. excursions. We exploit the regenerative property of $\{X_n\}$  to construct a non-parametric consistent estimator of the characteristic functions of $\rho_1$ and $\epsilon_1$ under the independence assumption  of $\rho_1$ and $\epsilon_1$.

%Let $\{X_n\}_{n\ge 0}$ be a sequence of random variables defined on a probability space and takes values in a measurable space $(S,\mathcal{S})$.
A sequence $\{X_n\}_{n\ge 0}$ is said to be delayed \textit{regenerative} if there exists a sequence $\{T_j\}_{j\geq 1}$ of positive integer valued random variables such that $\Pb(0<T_{j+1}-T_j<\infty)=1$ for all $j\geq 1$ and the random cycles $\eta_j\equiv (\{X_i: T_j\leq i< T_{j+1}\},T_{j+1}-T_j)$ for $j=1,2,\ldots$ are i.i.d. and independent of $\eta_0\equiv (\{X_i: 0\leq i< T_1\},T_{1})$. If $\{\eta_j\}_{j\geq 0}$ are i.i.d. then $\{X_n\}$ is called non-delayed regenerative sequence.
If, in addition, $E(T_2-T_1)<\infty$ then $\{X_n\}$ is called regenerative and positive recurrent.

If $\{X_n\}$ is a Markov chain with a general state space $(S,\mathcal{S})$, that is Harris irreducible and recurrent (see Definition 1) then it can be shown that $\{X_n\}$ is regenerative (\citet{AN1978}). Further if $\{X_n\}$ admits a stationary probability measure (necessarily unique because of irreducibility), then $\{X_n\}$ is  positive recurrent regenerative as well.

In Sections 2 and 3, under some condition on the distribution of $(\rho_1,\epsilon_1)$ we show that the sequence $\{X_n\}$  satisfying \eqref{model} is positive recurrent and regenerative by establishing that $\{X_n\}$ admits a stationary distribution and  is Harris irreducible, respectively. In Section 4, we show that the distribution of $(\rho_1,\epsilon_1)$ can be determined by transition probability function of $\{X_n\}$. We subsequently provide a consistent estimator of transition probability function of $\{X_n\}$ by using the regenerative property.
Finally, if $\rho_1$ and $\epsilon_1$ are independent then we provide a non-parametric consistent estimator of characteristic function of $\rho_1$ and $\epsilon_1$, based on $\{X_n\}_{n\geq 0}$.
%In this paper, we investigate the regenerative property of $X_n$ in \eqref{model} under two sets of conditions. In Section~3, Theorem~\ref{th7new} establishes the regeneration of $\{X_n\}$ when growth sequence $\{\rho_n\}$ has no mass at zero. A different approach is used to show regenerative property of $\{X_n\}$ in Theorem~\ref{theorem1} when $\Pb(\rho_1=0)>0$. In section~4, it is shown that the joint distribution of $(\rho_1,\epsilon_1)$ can be recovered from transition probability function of $X_n$ under some conditions.

\section{Limit distribution of $X_n$}
We begin with existence of the limiting distribution of $X_n$ in \eqref{model}.
\begin{theorem}\label{theorem2}
Let
%$\Pb(\rho_1=0)=0$
$-\infty\leq \E(\log|\rho_1|)<0$ and $\E(\log|\epsilon_1|)^{+}<\infty$. Then $\{X_n\}$ in \eqref{model} converges in distribution to $X_\infty$ as $n\to\infty$ where
\begin{equation}\label{eq2}
X_\infty\equiv\epsilon_1+\rho_1\epsilon_2+\rho_1\rho_2\epsilon_3+\ldots+\rho_1\ldots\rho_n\epsilon_{n+1}+\ldots.
\end{equation}
The infinite series on the right hand side of \eqref{eq2} is absolutely convergent with probability 1.
\end{theorem}
The above result can be deduced from \citet{B1986}. A proof of Theorem \ref{theorem2}
%under a slightly weaker hypothesis
is given in the appendix. Theorem~\ref{theorem2} does not indicate nature of limiting distribution of $X_n$. We show that the distribution of $X_\infty$ is non-atomic when the distribution of $(\rho_1,\epsilon_1)$ is non-degenerate.
\begin{theorem}\label{th5new}
Let $-\infty\leq \E(\log|\rho_1|)<0$, $\E(\log|\epsilon_1|)^{+}<\infty$, $\Pb(\rho_1=0)=0$ and $(\rho_1,\epsilon_1)$ has a non-degenerate distribution.
Then $X_\infty$ has a non atomic distribution, i.e., $\Pb(X_\infty=a)=0$ for all $a\in \mathbb R$.
\end{theorem}
\begin{proof} Since $(\rho_1,\epsilon_1)$ has a nondegenerate distribution,  the random variable  $X_\infty$ as in \eqref{eq2} does not have a degenerate distribution and hence  $\sup\{\Pb(X_\infty=a) : a\in\mathbb R\}\equiv p<1$. Let $a_0$ be such that $\Pb(X_\infty=a_0)=p$. Then, by Doob's martingale convergence theorem (see page 211 of \citet{AL2006}), we have
\begin{equation}\label{doobs result}\E(\mathbb I(X_\infty=a_0)|\mathcal F_n)\rightarrow \E( \mathbb I(X_\infty=a_0)|\mathcal F_\infty)
%=\mathbb I(X_\infty=a_0)
\quad  \mbox{w. p.}~ 1\end{equation}
where, $\mathcal F_n \equiv \sigma\{(\rho_i,\epsilon_i):i=1,2,\ldots,n, X_0\}$, the $\sigma$-algebra generated by $(\rho_i,\epsilon_i)$ for $i=1,\ldots, n$ and $X_0$, and $\mathcal F_\infty \equiv \sigma\{(\rho_i,\epsilon_i):i\in \mathbb N, X_0\}$.
Since $X_{\infty}$ is measurable with respect to $\mathcal F_\infty$, $\E( \mathbb I(X_\infty=a_0)|\mathcal F_\infty)=\mathbb I(X_\infty=a_0)$. Next \begin{eqnarray*}
&&\E(\mathbb I(X_\infty=a_0)|\mathcal F_n)\\
&=&\Pb(\epsilon_1+\rho_1\epsilon_2+\cdots+\rho_1\cdots\rho_{n-1}\epsilon_n+\rho_1\cdots \rho_n(\epsilon_{n+1}+\rho_{n+1}\epsilon_{n+2}+\cdots)=a_0|\mathcal F_n)\\
&=&\Pb\left(Y_n=\frac{a_0-\epsilon_1-\rho_1\epsilon_2-\cdots -\rho_1\rho_2\cdots \rho_{n-1}\epsilon_n}{\rho_1\rho_2\cdots \rho_n}|\mathcal F_n\right)\ \ (\mbox{since}\ \Pb(\rho_1=0)=0, \ |\rho_1\cdots \rho_n|\neq 0\ \forall\ n\geq 1)\\
\end{eqnarray*}
where $Y_n=\epsilon_{n+1}+\rho_{n+1}\epsilon_{n+2}+\rho_{n+1}\rho_{n+2}\epsilon_{n+3}+\cdots$. But $Y_n$ and $X_\infty$ have the same distribution, and $Y_n$ is independent of $\mathcal F_n$ and $\frac{a_0-\epsilon_1-\rho_1\epsilon_2-\cdots -\rho_1\rho_2\cdots \rho_{n-1}\epsilon_n}{\rho_1\rho_2\cdots \rho_n}$ is $\mathcal F_n$ measurable. So
$$\E(\mathbb I(X_\infty=a_0)|\mathcal F_n)\leq p<1 \ \mbox{for all}\ n\geq 1.$$
From \eqref{doobs result}, it follows that $\mathbb I(X_\infty=a_0)\leq p<1$ with probability 1. Since $\mathbb I(X_\infty=a_0)$ is a $\{0,1\}$ valued random variable, $\mathbb I(X_\infty=a_0)=0$ with probability 1 and hence $\Pb(X_\infty=a_0)=0$. Hence, $X_\infty$ has a non atomic distribution.
\end{proof}

A natural question is under what additional conditions on the distribution of $(\rho_1,\epsilon_1)$, the sequence $\{X_n\}$ is regenerative. When a Markov sequence is Harris recurrent and $\sigma$-algebra is countably generated then it can be established that the sequence exhibits regenerative property (see \citet{AN1978}). We now explore the Harris recurrence property of $\{X_n\}$.

\section{Harris recurrence of $X_n$}

\noindent\textbf{Definition~1.} A Markov chain $\{X_n\}_{n\geq 0}$ is called {\bf Harris or $\phi$-recurrent} if there exists a $\sigma$-finite measure $\phi$ on the state space $(S,\mathcal{S})$ such that
\begin{eqnarray}\label{def1}
\phi(A)>0 &\implies & \Pb(\tau_A<\infty|X_0=x)=1~\forall x\in S,
\end{eqnarray}
where $\tau_A=\min \{n:n\geq 1, X_n\in A\}$.

Note that any irreducible and recurrent Markov chain with a  countable state space is Harris recurrent as one can take $\phi$ to be the $\delta$ measure at some $i_0\in S$.
A definition related to Definition~1 is given by \citet{AN1978}.

\noindent \textbf{Definition~2.} A Markov chain $\{X_n\}$ is called $(A,\epsilon,\phi,n_0)$ \textbf{recurrent} if there exists a set $A\in\mathcal{S}$, a probability measure $\phi$ on $S$, a real number $\epsilon>0$, and an integer $n_0>0$ such that
\begin{equation}\label{def:cond1}
\Pb(\tau_A<\infty|X_0=x)\equiv \Pb_x(\tau_A<\infty)= \Pb_x(X_n\in A~\mbox{for some}~n\geq1)=1\quad\forall \  x\in S \ \mbox{and}\end{equation}
\begin{equation}\label{defn:connd2}
\Pb(X_{n_0}\in E|X_0=x)\equiv \Pb_x(X_{n_0}\in E)=  \Pb^{(n_0)}(x,E)\ge\epsilon\phi(E)\quad\forall~ x\in A~\mbox{and}~\forall \ E\subset \mathcal{S}.
\end{equation}
It can be shown by using the C-set lemma of Doob (see \citet{oreylecturenote}) that when $\mathcal{S}$ is countably generated, then Definition 1 implies Definition 2.
That Definition 2 implies Definition 1 is not difficult to prove.

The following theorem provides a sufficient condition for $\{X_n\}$ in \eqref{model} to be a Harris recurrent Markov chain.
\begin{theorem}\label{th7new}
%Let $ \E(\log|\rho_1|)<0$ and $\E(\log|\epsilon_1|)^{+}<\infty$.
%Let the hypothesis of Theorem \ref{theorem2} hold.
Let $-\infty\leq \E(\log|\rho_1|)<0$, $\E(\log|\epsilon_1|)^{+}<\infty$, $\Pb(\rho_1=0)=0$ and $-\infty<c<d<\infty$ be  such that $\Pb(c\le X_\infty\le d)>0$. Then, for all $ x\in \mathbb R$,
\begin{equation}\label{Harris:condition1}
\Pb(X_n\in [c,d]\ \mbox{ for some }\ n\geq 1|X_0=x)=1.
\end{equation}
In addition, let there exists a finite measure $\phi$ on $\mathcal{R}$ such that $\phi([c,d])>0$ and $0<\alpha<1$ such that
\begin{equation}\label{harris:condtion2}
\inf_{c\leq x\leq d}\Pb(\rho_1x+\epsilon_1\in \cdot)\geq \alpha \phi(\cdot).\end{equation}
Then,
the Markov chain $\{X_n\}_{n\geq 0}$ as described in (\ref{model}) is Harris recurrent and hence regenerative.
\end{theorem}
Note that since $X_n$ converges in distribution to  $X_{\infty}$ which is a proper real valued  random variable, $\{X_n\}_{n\geq 0}$  is positive recurrent as well. Thus under the hypothesis of Theorem \ref{th7new}, $\{X_n\}_{n\geq 0}$ is regenerative and positive recurrent.  The proof of Theorem \ref{th7new} is based on the following results.
\begin{lemma}\label{theorem4}
Let $-\infty\leq  \E(\log|\rho_1|)<0$, $\E(\log|\epsilon_1|)^{+}<\infty$, $\Pb(\rho_1=0)=0$ and  $-\infty<c<d<\infty$ be such that $\Pb(c\le X_\infty\le d)>0$.
Then, there exist $\theta>0$ and for all $x\in \mathbb R$, an integer $n_x\ge1$  such that
\begin{equation}\label{Harris irreducible}
\Pb(X_n\in[c,d]|X_0=x)\ge \theta \quad\mbox{for all}~ n\ge n_x.
\end{equation}
\end{lemma}
\begin{proof}
Iterating (\ref{model}) yields,
$$X_{n}=\rho_n \rho_{n-1}\cdots \rho_1 X_0+\rho_n \rho_{n-1}\cdots \rho_2 \epsilon_1+\cdots +\rho_n\epsilon_{n-1}+\epsilon_n\equiv Z_nX_0+Y_n,~~\mbox{say}.$$
So, if $X_0=x$ w. p. 1, then
\begin{eqnarray*}
\Pb_x(X_n\in[c,d])&=&\Pb(Y_n+Z_n x\in[c,d])\\
&\ge&\Pb(Y_n\in[c+\eta,d-\eta],|Z_n x|<\eta)\\
&\ge&\Pb(Y_n\in[c+\eta,d-\eta])-\Pb(|Z_n x|\ge\eta),
\end{eqnarray*}
where $\eta>0$ such that $c+\eta<d-\eta$. Now, define
\begin{equation}
Y'_n\equiv\epsilon_1+\rho_1\epsilon_2+\rho_1\rho_2\epsilon_3+\cdots+\rho_1\ldots\rho_{n-1}\epsilon_n.
\end{equation}
Note that the distribution of $Y_n$ and $Y'_n$ are same and from Theorem~\ref{theorem2}, $Y'_n\to X_\infty$ with probability 1. Thus, we have
\begin{eqnarray*}
\Pb_x(X_n\in[c,d])&\ge&\Pb(Y'_n\in[c+\eta,d-\eta])-\Pb(|Z_n x|\ge\eta)\\
&\ge&\Pb(Y'_n\in[c+\eta,d-\eta],|Y'_n-X_\infty|\le\eta^{'})-\Pb(|Z_n x|\ge\eta)\\
&\ge&\Pb(X_\infty\in[c+\eta+\eta^{'},d-\eta-\eta^{'}])-\Pb(|Y'_n-X_\infty|\ge\eta^{'})-\Pb(|Z_n x|\ge\eta),
\end{eqnarray*}
where $\eta^{'}>0$ such that $c+\eta+\eta^{'}<d-\eta-\eta^{'}$.

Now choose $n_1$ large such that  $\Pb(|Y'_n-X_\infty|\ge\eta^{'})\le\frac\delta2$ and $n_{2}$ large such that
$\Pb(|Z_{n_{2}} x|\ge\eta)\le\frac\delta2$. Note that choice of $n_2$ depends on $x$. Let $n_x=\max(n_1,n_{2})$. Then, for all $n\ge n_x$,
\begin{eqnarray*}
\Pb_x(X_{n_x}\in[c,d])\ge \Pb(X_\infty\in[c+\eta+\eta^{'},d-\eta-\eta^{'}])-\delta.
\end{eqnarray*}
Since $X_\infty$ has a continuous distribution by Theorem \ref{th5new} and $\Pb(c\le X_\infty\le d)>0$, first choose $\eta$ and $\eta^{'}$ and then $\delta$ small enough such that $$\theta\equiv\Pb(X_\infty\in[c+\eta+\eta^{'},d-\eta-\eta^{'}])-\delta>0.$$ Thus \eqref{Harris irreducible} is established.
%Now as $\eta, \eta'$ goes to zero, we have
%$$\Pb_x(X_{n_x}\in[c,d])\ge \Pb(X_\infty\in[c,d])-\delta,
%$$
%and choose $\delta$ small enough so that $\Pb(X_\infty\in[c,d])-\delta=\theta >0$.
\end{proof}

%We need the following lemma to prove Theorem~\ref{th7new}.
\begin{lemma}\label{lemma1}
Let $\{X_n\}$ be a time homogeneous Markov chain with state space $(S,\mathcal{S})$ and transition function $P(\cdot,\cdot)$. Let there exists $A\in \mathcal{S}$ and $0<\theta \leq 1$ such that for all $x\in S$, there exists an integer $n_x\geq 1$ such that
\begin{equation}\label{condition:lemma1}
\Pb(X_{n_x}\in A|X_0=x)\geq \theta .
\end{equation}
Then for all $x\in S$,
\begin{equation}\Pb(\tau_{A}<\infty|X_0=x)=1\end{equation}
where $\tau_{A}=\min \{n:\ n\geq 1, X_n\in A)\}$.
\end{lemma}
\begin{proof}
Fix $x\in S$. Let $B_0\equiv\{X_{n_x}\notin A\}$ and
$\tau_0=n_x$. Then $B_0\equiv \{X_{\tau_0}\notin A\}$.
Let us define
\begin{eqnarray*}
B_1&\equiv& \{X_{\tau_0}\notin A, \ X_{\tau_0+n_{X_{\tau_0}} }\notin A\}\\
\tau_1&=& \tau_0+n_{X_{\tau_0}}\\
B_2&\equiv & \{X_{\tau_0}\notin A,\ X_{\tau_1}\notin A,\ X_{\tau_1+n_{X_{\tau_1}} }\notin A \}\\
\tau_2&=& \tau_1+n_{X_{\tau_1}},
\end{eqnarray*}
and so on. Note $B_1=\{X_{\tau_0}\notin A, X_{\tau_1}\notin A\}$, $B_2=\{X_{\tau_0}\notin A, X_{\tau_1}\notin A, X_{\tau_2}\notin A\}$ and for any integer $k\geq 3$,
$$B_k\equiv \{ X_{\tau_0}\notin A, X_{\tau_1}\notin A,\ldots,X_{\tau_k}\notin A\},$$ with $\tau_k=\tau_{k-1}+n_{X_{\tau_{k-1}}}$. By hypothesis \eqref{condition:lemma1}, $\Pb((B_0)\leq (1-\theta)$. By the strong Markov property of $\{X_n\}$, $\Pb(B_1)\leq (1-\theta)^2$ and $\Pb(B_k) \leq (1-\theta)^{k+1}$ for all integer $k\geq 3$. This implies $\sum_{k=0}^\infty \Pb(B_k)<\infty$ since $\theta >0$. So $\sum_{k=0}^\infty\mathbb I_{B_k}(\cdot)<\infty$ with probability 1. This implies that with probability 1, $\mathbb I_{B_k}=0$ for all large $k>1$. That is, for all $x\in S$, $\Pb_x(X_{\tau_k}\in A\ \mbox{for some }k<\infty)=1$. Hence, for all $x\in S$, $\Pb_{x}(\tau_A<\infty)=1$.
\end{proof}

\begin{proof}[Proof of Theorem~{\ref{th7new}}]
In view of Definition 2, it is enough to prove \eqref{Harris:condition1} to show $\{X_n\}$ is Harris recurrent. The proof of \eqref{Harris:condition1} follows from Lemma~\ref{theorem4} and~\ref{lemma1}. Now from Lemma~2.2.5 of \citet{AA1998}, it follows that $X_n$ is regenerative.
\end{proof}

Theorem~{\ref{th7new}} provides sufficient conditions on $(\rho_1,\epsilon_1)$ so that the sequence $X_n$ becomes Harris recurrent and hence regenerative. These sufficient conditions are fairly general and hold for large class of distribution of $(\rho_1,\epsilon_1)$. Here are some examples where
\eqref{Harris:condition1} and \eqref{harris:condtion2} hold.

\noindent \textbf{Example 1:} $\epsilon_1$ is a standard normal, $N(0,1)$ random variable, $\rho_1$ has bounded support with $\E\log|\rho_1|<0$ and $\epsilon_1, \rho_1$ are independent.

\noindent \textbf{Example 2:} $\epsilon_1$ is a Uniform $(-1,1)$  random variable, $\rho_1$ has bounded support with $\E\log|\rho_1|<0$ and $\epsilon_1, \rho_1$ are independent.

In both the cases hypothesis of Theorem \ref{theorem2} hold and $X_\infty$ is of the form $(\epsilon_1+\rho_1 \tilde{X}_\infty)$ where $\tilde{X}_\infty$ has the same distribution as $X_\infty$ and independent of $X_\infty$. One can show in both above cases that for some $c<0<d$, $|c|$ and $d$ sufficiently small, conditions \eqref{Harris:condition1} and \eqref{harris:condtion2} hold.
%When growth sequence $\{\rho_n\}$ are independent of innovation sequence $\{\epsilon_n\}$, the following lemma provides a class of distribution function of $(\rho_1,\epsilon_1)$ for which conditions of Theorem~{\ref{th7new}} holds.

In Theorem~\ref{th7new}, growth sequence $\{\rho_n\}$ has no mass at zero and the regeneration property of $X_n$ is established by showing Harris recurrence of the sequence. When $\Pb(\rho_1=0)>0$, then the regenerative property of $X_n$ can be shown more easily.
\begin{theorem}\label{theorem1}
Let $\{X_n\}_{n\geq 0}$ be a RCAR(1) sequence as in \eqref{model}. If $\Pb(\rho_1=0)\equiv\alpha>0$, then $\{X_n\}_{n\geq 0}$ is a positive recurrent regenerative sequence.
\end{theorem}
\begin{proof}
Let $\tau_0=0$ and $\tau_{j+1}=\min\{n: n\ge \tau_j+1, \rho_n=0\}$ for $j\geq 0$. We need to show that
\begin{equation}\label{regeneration iid}
\Pb(\tau_{j+1}-\tau_j=k_j,X_{\tau_j+l}\in A_{l,j}, 0\leq l<k_j,1\leq j\leq r) =\prod_{j=1}^r \Pb(\tau_2-\tau_1=k_j,X_{\tau_1+l}\in A_{l,j},0\leq l<k_j)
\end{equation}
for all $k_1,k_2,\ldots,k_r\in \mathbb N$ and  $A_{l,j}\in \mathcal B(\mathbb R)$,
%the Borel sigma algebra of $\mathbb R$,
$0\leq l<k_j, j=1,2,\ldots,r$, $r=1,2,\ldots$.

Since $\{(\rho_n,\epsilon_n)\}_{n\geq 1}$ are i.i.d. and $\Pb(\rho_1=0)=\alpha>0$, it follows that $\{\tau_{j+1}-\tau_j,j\geq 0\}$ are i.i.d.  with jump distribution
$$\Pb(\tau_{j+1}-\tau_j=k)=(1-\alpha)^{k-1}\alpha, \quad \mbox{for}\quad k=1,2,\ldots,$$
that is, geometric with ``success" parameter $\alpha$. Next, since $\{(\rho_n,\epsilon_n)\}_{n\geq 1}$ are i.i.d. \eqref{regeneration iid} follows. Further since $\E(\tau_2-\tau_1)<\infty$, the sequence $\{X_n\}$ is positive recurrent regenerative.
\end{proof}

\noindent
{\bf Remark 1.} When $\Pb(\rho_1=0)=\alpha>0$ and the joint distribution of $(\rho_1,\epsilon_1)$ is discrete, then the limiting distribution $\pi$ of $X_\infty$ is a discrete probability distribution, that is, there exists a countable set $A_0$ in $\mathbb R^2$ such that $\pi(A_0)=1$. This is in contrast to Theorem \ref{th5new} which provides a sufficient condition for $X_\infty$ to have a non atomic distribution.

\section{Estimation of transition function and characteristic functions of $\rho_1$ and $\epsilon_1$}
The transition function $\Pb(x,A)$ of the Markov chain $\{X_n\}_{n\geq 0}$, defined by \eqref{model}, is precisely equal to $\Pb(\rho_1x+\epsilon_1\in A)$.  The following result determines the joint distribution of $(\rho_1,\epsilon_1)$ in terms of the transition function, $\Pb(\cdot,\cdot)$.
\begin{theorem}\label{theorem3}
If the distribution of $\rho_1x+\epsilon_1$ is known for all $x$ of the form $\frac{t_1}{t_2}$ where $t_2\neq 0$ and $(t_1,t_2)$ is dense in $\mathbb R^2$ then the distribution of $(\rho_1,\epsilon_1)$ is determined.
\end{theorem}
\begin{proof} For any $(t_1,t_2)\in \mathbb R^2$, the characteristic function of $(\rho_1,\epsilon_1)$ is
$$\psi_{(\rho_1,\epsilon_1)}(t_1,t_2)=\E(e^{i(t_1\rho_1+t_2\epsilon_1)})=\E(e^{it_2(\rho_1\frac{t_1}{t_2}+\epsilon_1)})=\phi_{\frac{t_1}{t_2}}(t_2)$$
where $\phi_x(t)=\E(e^{it(\rho_1x+\epsilon_1)})$ for all $x,t\in \mathbb R$. If $\phi_x(\cdot)$ is known for all $x$ of the form $\frac{t_1}{t_2}$ where $(t_1,t_2)$ is dense in $\mathbb R^2$, then $\psi_{(\rho_1,\epsilon_1)}(t_1,t_2)$ is determined for all such $(t_1,t_2)$ and hence by continuty for all $(t_1,t_2)\in \mathbb R^2$. Hence the distribution of $(\rho_1,\epsilon_1)$ is determined completely.
\end{proof}

Theorem \ref{theorem3} implies that if the transition function $\Pb(x,A)$ of $\{X_n\}_{n\geq 0}$ can be determined from observing the sequence sequence $\{X_n\}$, then the distribution of $(\rho_1,\epsilon_1)$ can also be determined.
%Also it suffices to know $\Pb(x,\cdot)$ for all $x$ of the form $\frac{t_1}{t_2},\ t_2\neq 0$ and $(t_1,t_2)$ dense in $\mathbb R^2$.
We now estimate the transition probability function $\Pb(x,(-\infty,y])$, for $x,y\in \mathbb R^2$ from the data $\{X_i\}_{i=0}^n$. In the following theorems, we show that the estimator $F_{n,h}(x,y)$, given in \eqref{eqn:estimator of transition fn} below, is a strongly consistent estimator for $\Pb(X_1\leq y|X_0=x)$.

\begin{theorem}\label{estimator of transition fn}
Let $\{X_n\}$ satisfies the hypothesis of Theorem \ref{th7new}. For $n\geq 1$, $h>0,\ x,y\in \mathbb R$, let
\begin{equation}\label{eqn:estimator of transition fn}
F_{n,h}(x,y)=
\begin{cases}
\frac{\frac{1}{nh}\sum_{i=0}^{n-1}\mathbb I(x\leq X_i \leq x+h, X_{i+1}\leq y)}{\frac{1}{nh}
\sum_{i=0}^n \mathbb I(x\leq X_i\leq x+h)} & \mbox{ if } \ \ \mathbb I(x\leq X_i\leq x+h)\neq 0 \mbox{ for some }0\leq i \leq n-1\\
0 & \mbox{ otherwise,}
\end{cases}
\end{equation}
where $\mathbb I(A)$ denotes the indicator function of the event $A$.
\begin{enumerate}[(a)]
\item Then with probability $1$, for each $x,y\in \mathbb R$
\begin{equation}\label{lim of F as n}
\lim_{n\to \infty}  F_{n,h}(x,y) \equiv \psi(x,y,h)= \frac{\int_x^{x+h}G(u,y) \Pb(X_\infty \in du )}{\Pb(X_\infty \in (x,x+h])},
\end{equation}
where $G(u,y)= \Pb(x, (-\infty,y])=\Pb(X_1\leq y|X_0=x)$.
\item In addition, let $G(x,y)$ and the random variable $X_\infty$ satisfy
\begin{equation}\label{sufficient condition limF_n,h as h}
\lim_{h\to 0} \frac{\int_x^{x+h} G(u,y) \Pb(X_\infty \in du)}{\Pb(x<X_\infty\leq x+h)}=G(x,y), \ \mbox{ for } \ x,y\in \mathbb R^2.
\end{equation}
Then for $x,y\in \mathbb R$
\begin{equation}\label{lim of F as n,h}
\lim_{h\to 0}\lim_{n\to 0} F_{n,h}(x,y)=\Pb(X_1\leq y|X_0=x), \mbox{ with probability }1.
\end{equation}
 \end{enumerate}
\end{theorem}
\begin{proof}
Since $\{X_i\}_{i\ge0}$ is regenerative and positive recurrent, the vector sequence $\{(X_i, X_{i+1})\}_{i\ge0}$ is also regenerative and positive recurrent Markov chain. The numerator in \eqref{eqn:estimator of transition fn} converges to $\int_x^{x+h}G(u,y) \Pb(X_\infty \in du)$ with probability 1 by using Theorem~9.2.10 of \citet{AL2006}. Similarly denominator converges to $\Pb(X_\infty \in (x,x+h])$ with probability 1. This completes the proof of part (a).

The proof of part (b) follows from \eqref{lim of F as n} and \eqref{sufficient condition limF_n,h as h}.
\end{proof}
\noindent
\textbf{Remark 2.} A sufficient condition for \eqref{sufficient condition limF_n,h as h} to hold  is that  the distribution of $X_\infty$ is absolutely continuous with strictly positive and continuous density function and  the function $G(x,y)$ is continuous in $x$ for fixed $y$.

%Now if we let $h\downarrow 0$ the right hand side of \eqref{eqn:estimator of transition fn} will converge to the transition probability function.
%\begin{theorem}
%Suppose $\{X_n\}$ satisfies the hypothesis of Theorem \ref{th7new} and in addition if $X_\infty$ has a probability density function $f_\infty$ such that $f_\infty(x)>0$,
%\end{theorem}
%Proof of the above theorem is straightforward, we omit it.

The following result is similar to that of Theorem \ref{estimator of transition fn}.
\begin{theorem}\label{estimation of ch fn}
Fix $x,t,h\in \mathbb R$. Let
\begin{equation*}
\phi_{n,h,x}(t)=
\begin{cases}
\frac{\frac{1}{nh}\sum_{j=0}^{n-1}e^{itX_{j+1}}\mathbb I(x <  X_j \leq x+h)}{\frac{1}{nh}
\sum_{j=0}^{n-1} \mathbb I(x <  X_j\leq x+h)} &\ \mbox{ if } \ \ \mathbb I(x\leq X_i\leq x+h)\neq 0 \mbox{ for some }0\leq i \leq n-1,\\
0 & \mbox{ otherwise.}
\end{cases}
\end{equation*}
Then
\begin{equation*}
\lim_{h\to 0}\lim_{n\to \infty} \phi_{n,h,x}(t)=\E(e^{it(\rho_1 x+\epsilon_1)}) \ \mbox{with probability } 1,
\end{equation*}
provided
\begin{equation}\label{phi limit as h}
\lim_{h\to 0}\frac{ \frac{1}{h}\int_x^{x+h} \E(e^{it(\rho_1 x+\epsilon_1)})\Pb(X_\infty \in du)}{\frac{1}{h}\int_x^{x+h}\Pb(X_\infty \in du)}=\E(e^{it(\rho_1 x+\epsilon_1)}).
\end{equation}
\end{theorem}
\begin{proof}Proof of this theorem is similar to the proof of Theorem \ref{estimator of transition fn} and hence omitted.\end{proof}

\noindent
\textbf{Remark 3.} A sufficient condition for \eqref{phi limit as h} to hold  is that  the distribution of $X_\infty$ is absolutely continuous with strictly positive and continuous density function on $(-\infty,\infty)$ and  the function $\E(e^{it(\rho_1 x+\epsilon_1)})$ is continuous in $x$ for fixed $t$.

Let  $\{\rho_1\}$ and $\{\epsilon_1\}$ are independent random variables.Then
%\begin{theorem}
%Under the Assumptions,.....
%\begin{eqnarray*}
%\lim_{h\to 0}\lim_{n\to \infty}\phi_{n,h,x}(t) &=&\E_x(e^{itX_1}), \ \mbox{for all} \ x,t\in \mathbb R \ \ \mbox{with probability}~1.
%\end{eqnarray*}
%\end{theorem}
\begin{equation*}
\phi_x(t)\equiv \E_x(e^{itX_1})= \E(e^{it(\rho_1x+\epsilon_1)})=\psi_{\rho}(tx) \psi_{\epsilon}(t),
\end{equation*}
where $\psi_{\rho}(t)=\E(e^{it\rho})$ and $\psi_{\epsilon}(t)=\E(e^{it\epsilon})$.
Also, note that
\begin{equation*}
  \psi_{\epsilon}(t) = \phi_0(t) \ \ \mbox{and}\ \
  \psi_{\rho}(tx)= \frac{\phi_x(t)}{\phi_0(t)}, \ \mbox{when} \ \ \psi_{\epsilon}(t)\neq0.
\end{equation*}
This yields the following corollary of Theorem \ref{estimation of ch fn}.
 \begin{corollary}
Let $\rho_1$ and $\epsilon_1$ be independent and conditions of Theorem \ref{estimation of ch fn} holds. Then
\begin{enumerate}[(a)]
\item $
\lim_{h\to 0}\lim_{n\to \infty}\phi_{n,h,0}(t)=\psi_{\epsilon}(t) \ \mbox{for all} \ t\in \mathbb R \ \ \mbox{with probability}~1.$
\item Let $\psi_{\epsilon}(t)\neq0$ for all $t\in \mathbb R$, then for all $x\neq0$
\begin{align*}
\lim_{h\to 0}\lim_{n\to \infty}\frac{\psi_{n,h,x}(t/x)}{\phi_{n,h,0}(t/x)}&=\psi_{\rho}(t) \ \mbox{for all} \ t\in \mathbb R \ \ \mbox{with probability}~1.
\end{align*}
\end{enumerate}
\end{corollary}

\bigskip\noindent
\textbf{Acknowledgement.} K. B. Athreya would like to thank the Department of Mathematics, IIT Bombay and in particular, Prof. Sudhir Ghorpade for offering him visiting professorship.
% during July 1 to August 31,  2014 and July 1 to August 31, 2015 where the work presented here was done.
 The work of Radhendushka Srivastava is partially supported by INSPIRE research grant of DST, Govt. of India and a seed grant from IIT Bombay.

\section{Appendix}
\begin{proof}[Proof of Theorem \ref{theorem2}:] Choose $\epsilon>0$ such that $\E(\log|\rho_1|)+\epsilon<0$. Now, by the strong law of large number,
$$\E(\log|\rho_1|)<0 \Rightarrow \frac{1}{n}\sum_{i=1}^n \log|\rho_i|\leq \E(\log |\rho_1|)+\epsilon, $$
 for sufficiently large $n$, with probability 1. Hence
\begin{equation}\label{eqn:bound prod rho}
|\rho_1\rho_2\ldots \rho_n|\leq e^{-n\lambda},
\end{equation}
 where $0<\lambda\equiv -(\E(\log|\rho_1|+\epsilon)<\infty$, for all large $n$, with probability 1.

Also $\E(\log|\epsilon_1|)^{+}<\infty$ implies that for any $\mu>0$, $\sum_{n=1}^\infty \Pb(\log|\epsilon_1|>n\mu )<\infty$ and hence  $\sum_{n}\Pb(\log|\epsilon_n|>n\mu)<\infty$. By Borel Cantelli lemma,  $|\epsilon_n|\leq e^{n\mu}$ for all $n$ large enough, with probability 1.

Now choose $0<\mu<\lambda$. Then, for sufficiently large $n$, with probability 1,
$$|\epsilon_{n+1} \rho_1\rho_2\ldots \rho_n|\leq e^{-n\lambda}e^{(n+1)\mu}.$$
Therefore $\sum_n |\epsilon_{n+1}|\rho_1\rho_2\ldots \rho_n|<\infty$ with probability 1.   Hence $\tilde{X}_\infty=\epsilon_1+\rho_1\epsilon_2+\rho_1\rho_2\epsilon_3+\ldots+\rho_1\ldots\rho_n\epsilon_{n+1}+\ldots
$ is well defined.

 Observe that
\begin{eqnarray*}
X_{n}&=&\rho_{n}(\rho_{n-1}X_{n-2}+\epsilon_{n-1})+\epsilon_n\\
&=&\rho_n \rho_{n-1}\cdots \rho_1 X_0+\rho_n \rho_{n-1}\cdots \rho_2 \epsilon_1+\cdots +\rho_n\epsilon_{n-1}+\epsilon_n
\end{eqnarray*}
and which has the same distribution as
\begin{equation}\label{eqn:init_term}
\epsilon_1+\rho_1\epsilon_2+\cdots+\rho_1\rho_2\cdots\rho_{n-1}\epsilon_n+\rho_1 \rho_{2}\cdots \rho_{n} X_0.
\end{equation}
Now by using (\ref{eqn:bound prod rho}) and above, we have $|\rho_1 \rho_{2}\cdots \rho_{n} X_0|$ converges to zero with probability 1. Thus, from (\ref{eqn:init_term}), as $n\to\infty$, we have
$$X_{n}\stackrel{ d}{\to}\tilde{X}_\infty,$$
where $\stackrel{ d}{\to}$ stands for  convergence in distribution.
\end{proof}

\bibliographystyle{abbrvnat}
%\bibliography{reference}

\end{document}